\documentclass[11pt,a4paper]{amsart}
\usepackage{graphicx}
\usepackage{amssymb}
\usepackage{amsmath}
\usepackage{amsthm}
\usepackage{amscd}
\usepackage[all,2cell,dvips]{xy}
\UseAllTwocells \SilentMatrices
\newtheorem{thm}{Theorem}[section]

\newtheorem{cor}[thm]{Corollary}
\newtheorem{lem}[thm]{Lemma}

\newtheorem{prop}[thm]{Proposition}
\theoremstyle{definition}
\newtheorem{defn}[thm]{Definition}
\theoremstyle{remark}

\numberwithin{equation}{section}

\begin{document}
\title[Noncommutative Gorenstein Projective Schemes]
{Noncommutative Gorenstein Projective Schemes and
Gorenstein-Injective Sheaves}
\author[  Xiao-Wu Chen
] {Xiao-Wu Chen}
\thanks{This project was supported by China Postdoctoral Science Foundation No. 20070420125, and
was also partly supported by the National Natural Science Foundation
of China (Grant No.s 10501041 and 10601052).}
\thanks{E-mail:
xwchen$\symbol{64}$mail.ustc.edu.cn}
\maketitle
\date{}%
\dedicatory{}%
\commby{}%
\begin{center}
Department of Mathematics\\
 University of Science and
Technology of China \\Hefei 230026, P. R. China
\end{center}

\begin{abstract}
We prove that if a positively-graded ring $R$ is Gorenstein and the
associated torsion functor has finite cohomological dimension, then
the corresponding noncommutative projective scheme ${\rm Tails}(R)$
is a Gorenstein category in the sense of \cite{EEG}. Moreover, under
this condition, a (right) recollement relating Gorenstein-injective
sheaves in ${\rm Tails}(R)$ and  (graded) Gorenstein-injective
$R$-modules is given.
\end{abstract}

\section{Introduction}

\subsection{} Gorenstein-projective and Gorenstein-injective modules
over arbitrary rings are introduced to generalize Auslander's notion
of modules of G dimension zero, and later they play central roles in
the theory of relative homological algebra  and Tate cohomology (see
\cite{Buc, EJ, Be2}). If the ring is Gorenstein, or more generally,
the module category is a Gorenstein category, then these modules
behave nicely. Note that via old and new works, those modules appear
naturally in the study of singularity categories of rings and
schemes (\cite{Buc, Ha2, Be2, Kr, CZ, C}).

\vskip 5pt

 Extensions of these notions to Grothendieck categories are
considered in \cite{EEG}. Since a Grothendieck category has enough
injective objects but usually not enough projective objects,
 in these categories we have the notion of Gorenstein-injective objects but
perhaps not of Gorenstein-projective objects. The notion of a
Gorenstein category is defined in \cite{EEG}. This is a Grothendieck
category with a generator of finite projective dimension, and
objects of finite projective dimension coincide with those of finite
injective dimension, and where these dimensions are uniformly
bounded. Again, when these Gorenstein condition is fulfilled, the
Gorenstein-injective objects behave nicely (see \cite[Theorem
2.24]{EEG} and our Proposition \ref{keryproposition}). Two main
types of examples of Gorenstien categories are the module category
$R\mbox{-Mod}$ over Gorenstein rings $R$ and the category  of
quasi-coherent sheaves ${\rm Qcoh}(\mathbb{X})$ for a locally
Gorenstein closed subscheme $\mathbb{X}\subseteq \mathbb{P}^n(k)$
for any commutative noetherian ring $k$ (see \cite{EEG}). In this
paper, we will show that the category of quasi-coherent sheaves on
certain noncommutative projective scheme in the sense of  \cite{AZ}
is Gorenstein, and under this Gorenstein condition, a (right)
recollement on the stable category of Gorenstein-injective sheaves
will be given explicitly.

\subsection{} We will state our main results in this subsection. Let
$R=\oplus_{n\geq 0}R_n$ be a positively-graded ring. Assume that $R$
is left-noetherian, and note that this is equivalent to that $R$ is
graded left-noetherian (i.e., any graded left ideal is
finitely-generated). Let $R\mbox{-Gr}$ be the category of left
$\mathbb{Z}$-graded $R$-modules with morphisms preserving degrees.
Recall the degree-shift functors $(d): R\mbox{-Gr} \longrightarrow
R\mbox{-Gr}$ for each integer $d$: let $M=\oplus_{n\in \mathbb{Z}}
M_n$ be a graded module, define $M(d)=\oplus_{n\in
\mathbb{Z}}M(d)_n$ such that $M(d)_n=M_{n+d}$, and the $R$-module
structure is unchanged; $(d)$ acts on morphisms naturally. Then
these functors are automorphisms.  Let $M=\oplus_{n\in \mathbb{Z}}
M_n$ be a graded $R$-module, $m\in M$. We say that $m$ is a
\emph{torsion element} if there exists $d\geq 1$ such that $R_{\geq
d}.m=0$, where $R_{\geq d}:=\oplus_{n\geq d}R_n$. Denote by
$\tau(M)\subseteq M$ the submodule consisting of all the torsion
elements, and it is a graded submodule. A module $M$ is said to be a
\emph{torsion module} if $\tau(M)=M$. Denote by $R\mbox{-Tor}$ the
full subcategory of $R\mbox{-Gr}$ consisting of torsion modules. It
is a Serre subcategory, i.e., it is closed under (graded)
submodules, quotient modules and extensions. The operation $\tau$
induces a functor $\tau: R\mbox{-Gr} \longrightarrow R\mbox{-Tor}$,
called the \emph{torsion functor},  which is the right adjoint to
the inclusion functor ${\rm inc}: R\mbox{-Tor} \longrightarrow
R\mbox{-Gr}$, and in other words, $R\mbox{-Tor}$ is a localizing
Serre subcategory. By \cite[Proposition 2.2(2)]{AZ}, the subcategory
$R\mbox{-Tor}$ is closed under essential monomorphisms, that is,
given an essential monomorphism $M \longrightarrow M'$ in
$R\mbox{-Gr}$, if $M$ is a torsion module so is $M'$. In summary,
our Setup 4.2 in section 4 applies in this situation.

\vskip 5pt

 Use the above notation. Artin and Zhang define the following quotient
 abelian category in the sense of Gabriel
  \begin{align*}
{\rm Tails}(R) \; :=\; R\mbox{-Gr}/{R\mbox{-Tor}}.
  \end{align*}
We denote by $\pi: R\mbox{-Gr}\longrightarrow {\rm Tails}(R)$ the
canonical functor, and set $\mathcal{O}=\pi(R)$ (see \cite{AZ} for
details). In this case, $\pi$ admits a right adjoint which will be
denoted by $\omega: {\rm Tails}(R) \longrightarrow R\mbox{-Gr}$.
Since the subcategory $R\mbox{-Tor}$ is closed under the
degree-shift functors, thus we get the induced functors $(d): {\rm
Tails}(R) \longrightarrow {\rm Tails}(R)$ for each $d$, which are
still automorphisms. The category ${\rm Tails}(R)$ is sometimes
called the \emph{ category of tails} of $R$, or it is known as the
noncommutative projective scheme: more precisely, write
$\mathbb{X}=({\rm Tails(R)}, \mathcal{O}, (1))$, then we call the
triple $\mathbb{X}$ the \emph{(polarized) general projective
scheme}, and $\mathcal{O}$ the \emph{structure sheaf}, $(1)$ the
\emph{twist functor}, we even write ${\rm Qcoh}(\mathbb{X})={\rm
Tails}(R)$, which is referred as the category of
\emph{quasi-coherent sheaves} on $\mathbb{X}$, and the graded ring
$R$ is called the \emph{homogeneous coordinate ring} of
$\mathbb{X}$. \vskip 5pt

We will be concerned with the question of when the category ${\rm
Tails}(R)$ is Gorenstein. A natural condition is that the ring $R$
is Gorenstein: recall that a ring $R$ is said to be
\emph{Gorenstein} if $R$ is two-sided noetherian and $R$ has finite
injective dimension both as the left and right regular module.
Similarly one defines the notion of \emph{graded Gorenstein rings}.
Thanks to a remark by Van den Bergh (\cite[p.670, line 9-11]{Van}),
we infer that a graded module $M$ has finite injective dimension in
$R\mbox{-Gr}$ if and only if $M$ has finite injective dimension as
an ungraded module. Thus the graded ring $R$ is graded Gorenstein if
and only if it is Gorenstein as an ungraded ring; and in this case,
one sees easily that $R\mbox{-Gr}$ is a Gorenstein category (by
comparing or applying the results in \cite[section 9.1]{EJ}).
Another condition we want to impose is that the functor $\tau$ has
finite cohomological dimension, that is, there exists $d\geq 1$ such
that the $n$-th right derived functor $R^n\tau=0$ for $n\geq d$. The
condition is satisfied in  commutative algebraic geometry, more
precisely, if the graded noetherian ring $R$ is commutative,
generated by $R_0$ and $R_1$, and $R_0$ is an affine algebra over
some field, then the functor $\tau$ always has finite cohomological
dimension (to see this well-known fact, note that $R$ has finite
Krull dimension, first apply the famous Serre theorem (\cite[p.125,
Ex. 5.9]{Har2}) and then Grothendieck's vanishing theorem
(\cite[p.208]{Har2}) to obtain that the structure sheaf
$\mathcal{O}$ has finite projective dimension, and then use our
Lemma \ref{cohomologicaldimensioncondition}). More generally, for
any noetherian commutative graded ring $R$ the torsion functor
$\tau$ has finite cohomological dimension. To see this, first note
that $\tau$ coincides with the local cohomology functor with respect
to the ideal $R_{\geq 1}$ (compare \cite[section 4]{Van}) and then
we can compute local cohomology via the Koszul complex of a finite
set of generators of the ideal $R_{\geq 1}$ (\cite[Theorem 2.3 and
Theorem 2.8]{Gro}).

\vskip 5pt

 We have our first result, which is a direct consequence of Corollary
 \ref{improtantcorollary}. One may compare it with the main example  of \cite[Theorem
 3.8]{EEG} in commutative algebraic geometry.

\vskip 10pt

 \noindent {\bf Theorem A}\quad \emph{Let $R=\oplus_{n\geq 0}R_n$ be a graded
Gorenstein ring such that the associated torsion functor $\tau$ has
finite cohomological dimension. Then the category of tails ${\rm
Tails}(R)$ is a Gorenstein category. }

\vskip 10pt

Next we will be concerned with the full subcategory of
Gorenstein-injective objects (or sheaves) in ${\rm Tails}(R)$ and
more precisely its stable category modulo injective objects. Recall
that in an abelian category $\mathcal{A}$ with enough injective
objects, an object $G$ is said to be Gorenstein-injective
(\cite{EJ}), if there exists an exact complex $I^\bullet$ of
injective objects such that the Hom complex ${\rm
Hom}_\mathcal{A}(Q, I^\bullet)$ is still exact for any injective
object $Q$, and $G=Z^0(I^\bullet)$ is the $0$-th cocycle. The full
subcategory of Gorenstein-injective objects will be denoted by ${\rm
GInj}(\mathcal{A})$. Note that ${\rm GInj}(\mathcal{A})$ is a
Frobenius exact category and its relative injective-projective
objects are precisely the injective objects in $\mathcal{A}$ (see
our Lemma \ref{basicproperty}). Then the stable category
$\underline{{\rm GInj}}(\mathcal{A})$ modulo injective objects has a
canonical triangulated structure by \cite[Chapter 1, section
2]{Ha1}. \vskip 5pt

We have our second result which relates the category of
Gorenstein-injective objects (or sheaves) in ${\rm Tails}(R)$ to the
category of Gorenstein-injective $R$-modules.  It is a direct
consequence of Theorem \ref{maintheoremII}. \vskip 10pt

 \noindent {\bf Theorem B}\quad \emph{Use
the same assumption in Theorem A. Then we have a right recollement
of triangulated categories:}
\[\xymatrix@C=40pt{
{\underline{\rm GInj}(R\mbox{\rm -Tor})} \ar@<+.7ex>[r]^-{\rm inc} &
\ar@<+.7ex>[l]^-{\underline{\tau}} {\underline{\rm GInj}(R\mbox{\rm
-Gr})} \ar@<+.7ex>[r]^-{\underline{\pi}} &
\ar@<+.7ex>[l]^-{\underline{\omega}} \underline{\rm GInj}({\rm
Tails}(R)), }\] \emph{where the functors are induced from the ones
between the abelian categories.}

\vskip 10pt

Let us remark that we actually obtain a recollement in this
situation. By Theorem A and Proposition \ref{keryproposition}, we
observe that an acyclic (= exact) complex of injective objects in
${\rm Tails}(R)$ is totally-acyclic, and then by \cite[Proposition
7.2]{Kr}, we infer that $\underline{\rm GInj}({\rm Tails}(R))$ is
triangle equivalent to the \emph{stable derived category} $S({\rm
Tails}(R))$ in the sense of Krause (\cite[section 5]{Kr}). In
particular, by noting that the Grothendieck category ${\rm
Tails}(R)$ is locally noetherian and then by  \cite[Corollary
5.4]{Kr}, the stable category $\underline{\rm GInj}({\rm Tails}(R))$
is compactly generated. Recall a well-known natural isomorphism (see
\cite[equation (3.12.2)]{AZ}), for any sheaf $\mathcal{M}\in {\rm
Tails}(R)$, $\omega(\mathcal{M})\simeq \oplus_{d\in \mathbb{Z}}{\rm
Hom}_{{\rm Tails}(R)}(\mathcal{O}, \mathcal{M}(d))$, and since the
structure sheaf $\mathcal{O}$ is a noetherian object, we infer that
the functor $\omega$ preserves arbitrary coproducts, and
consequently, the induced triangle functor $\underline{\omega}:
\underline{\rm GInj}({\rm Tails}(R)) \longrightarrow \underline{\rm
GInj}(R\mbox{-Gr})$ also preserves arbitrary coproducts. Applying
Brown representability theorem (\cite[Theorem 8.4.4]{Ne}) to the
functor  $\underline{\omega}$, we obtain a right adjoint functor
$\pi'$ for $\underline{\omega}$, and thus consequently, we obtain
also a right adjoint functor $i'$ for $\underline{\tau}$, and both
$\pi'$ and $i'$ are triangle functors (see the arguments in {\bf
2.2}).  In summary, we obtain a recollement in this situation,
expressed in the following diagram

\[\xymatrix@C=40pt{
\underline{\rm GInj}({\rm Tails}(R))
\ar[r]^-(.45){\underline{\omega}} &
\ar@<-2.0ex>[l]_{\underline{\pi}} \ar@<+2.0ex>[l]_(.48){\pi'}
{\underline{\rm GInj}(R\mbox{\rm -Gr})} \ar[r]^-{\underline{\tau}} &
\ar@<-2.0ex>[l]_{\rm inc} \ar@<+2.0ex>[l]_(.49){i'}  {\underline{\rm
GInj}(R\mbox{\rm -Tor})}.
 }\]

\vskip 15pt

\subsection{} The paper is organized as follows : in section 2, we
review Gabriel's theory of quotient abelian categories and Verdier's
theory of quotient triangulated categories, and their right
recollements; in section 3, we recall the definition of Gorenstein
categories in \cite{EEG} with slight modifications and we
characterize the Gorenstein condition by the bounded derived
category; in section 4, we consider when the quotient abelian
category of a Gorenstein category is still Gorenstein, and we obtain
our main results Theorem \ref{maintheorem} and Corollary
\ref{improtantcorollary}, which give our Theorem A; in the final
section, we consider the stable category of Gorenstein-injective
objects in abelian categories, and we obtain a right recollement
relating the stable categories of Gorenstein-injective objects in
the abelian category and its quotient category, see Theorem
\ref{maintheoremII}, which gives our Theorem B.

\section{Preliminaries}

\subsection{} Let us recall Gabriel's  theory of quotient abelian
categories. Let $\mathcal{A}$ be an arbitrary abelian category. Let
$\mathcal{B}\subseteq \mathcal{A}$ be a \emph{Serre subcategory},
that is, $\mathcal{B}\subseteq \mathcal{A}$ is a full subcategory
which is closed under subobjects, quotient objects and extensions,
or equivalently, for any short exact sequence $0\longrightarrow X
\longrightarrow Y \longrightarrow Z \longrightarrow 0$ in
$\mathcal{A}$, $Y$ lies in $\mathcal{B}$ if and only if $X$ and $Z$
do. Recall the \emph{quotient category} $\mathcal{A}/\mathcal{B}$ is
defined by taking the inverse, formally, of the morphisms $f:
X\longrightarrow Y $ in $\mathcal{A}$ such that both ${\rm Ker}\; f$
and ${\rm Coker}\; f$ lie in $\mathcal{B}$. More precisely, denote
by $\Sigma$ the class of such morphisms, which is a multiplicative
system on $\mathcal{A}$, then the category $\mathcal{A}/\mathcal{B}$
is given as follows: its objects coincide with the ones in
$\mathcal{A}$; and morphisms are given by (right) fractions (=
roofs) $f/s: X \stackrel{s}\Longleftarrow Z
\stackrel{f}\longrightarrow Y$ modulo certain equivalence relation,
where $s\in \Sigma$; composition of fractions are given in
\cite[Appendix A.2]{Ne}.  Then we have a canonical functor $\pi:
\mathcal{A}\longrightarrow \mathcal{A}/\mathcal{B}$. Note that $\pi$
sends morphism $f: X\longrightarrow Y$ to the trivial fraction
$f/{{\rm Id}_X}: X \stackrel{{\rm Id}_X}\Longleftarrow X
\stackrel{f}\longrightarrow Y$. Let us remark that the morphism
space in $\mathcal{A}/\mathcal{B}$ can also be expressed as certain
colimits of morphism spaces in $\mathcal{A}$ as follows
\begin{align*}
{\rm Hom}_{\mathcal{A}/\mathcal{B}}(\pi(X), \pi(Y))\; = \; {\rm
colimt}_{X',Y'} \; {\rm Hom}_\mathcal{A}(X', Y/Y'),
\end{align*}
where $X'$ and $Y'$ run over all the subobjects of $X$ and $Y$ such
that $X/X'\in \mathcal{B}$ and $Y'\in \mathcal{B}$.
 \vskip 5pt

We collect some needed properties of the quotient category: the
category $\mathcal{A}/\mathcal{B}$ is abelian and $\pi$ is an exact
functor with ${\rm Ker}\; \pi=\mathcal{B}$, or in other words, there
is a   short ``exact sequence" of abelian categories:
\begin{align*}
\mathcal{B} \stackrel{\rm inc}\longrightarrow \mathcal{A}
\stackrel{\pi}\longrightarrow \mathcal{A}/\mathcal{B}
\end{align*}
where ${\rm inc}$ is the inclusion functor; moreover, any short
exact sequence in $\mathcal{A}/\mathcal{B}$ is isomorphic to the
image of some short exact sequence in $\mathcal{A}$.

\vskip 5pt

Later we will be mainly interested in the localizing Serre
subcategories. Recall that a Serre subcategory $\mathcal{B}\subseteq
\mathcal{A}$ is said to be \emph{localizing} provided that the
inclusion functor ${\rm inc}: \mathcal{B} \hookrightarrow
\mathcal{A}$ admits a right adjoint. This adjoint functor will be
denoted by $\tau$ and called the \emph{torsion functor}. In this
case, the quotient functor $\pi: \mathcal{A}\longrightarrow
\mathcal{A}/\mathcal{B}$ also admits a right adjoint, denoted by
$\omega: \mathcal{A}/{\mathcal{B}}\longrightarrow \mathcal{A}$.
Observe that both $\tau$ and $\omega$ are fully-faithful functors.
Consider the full subcategory $\mathcal{B}^\perp := \{X\in
\mathcal{A}\; |\; {\rm Hom}_\mathcal{A}(\mathcal{B}, X)={\rm
Ext}^1_{\mathcal{A}}(\mathcal{B}, X)=0\}$, whose objects are said to
be $\mathcal{B}$-\emph{local} (\cite[Appendix A.2]{Ne}). In this
case $\mathcal{B}^\perp$ coincides with the essential image ${ \rm
Im}\; \omega$ of the functor $\omega$. A remarkable fact is that the
full subcategory $\mathcal{B}^\perp$ is a Giraud subcategory of
$\mathcal{A}$. Recall that a full subcategory $\mathcal{B}'$ is said
to be a \emph{Giraud subcategory} \cite{St} if the inclusion functor
$\mathcal{B'}\hookrightarrow \mathcal{A}$ admits a left adjoint
functor $\sigma: \mathcal{A} \longrightarrow \mathcal{B}'$ which
preserves kernels. Conversely, given a Giraud subcategory
$\mathcal{B}'$, the kernel ${\rm Ker}\; \sigma$ which coincides with
$^\perp\mathcal{B'}:=\{X\in \mathcal{B}\; |\; {\rm
Hom}_\mathcal{A}(X, \mathcal{B}')=0\}$ is a localizing Serre
subcategory. In fact, this gives a bijection between the class of
Serre subcategories and the class of Giraud subcategories.

 \vskip
5pt

Another related concept is the right recollement of abelian
categories. Let $\mathcal{A},\mathcal{A}', \mathcal{A}''$ be abelian
categories. A \emph{right recollement} of $\mathcal{A}$ with respect
to $\mathcal{A}'$ and $\mathcal{A}''$ is expressed by the following
diagram of functors
\[\xymatrix@C=40pt{
\mathcal{A}' \ar@<+.7ex>[r]^-{i_{!}} & \ar@<+.7ex>[l]^-{i^{!}}
\mathcal{A} \ar@<+.7ex>[r]^-{j^*} & \mathcal{A}''
\ar@<+.7ex>[l]^-{j_*} }\] satisfying the following conditions\\
(1).\quad The pairs $(i_!, i^!)$ and $(j^*, j_*)$ are adjoint
pairs.\\
(2).\quad The functors $i_!$ and $j^*$ are exact such that $j^*i_!=0$. \\
(3).\quad The functors $i_!$ and $j_*$ are full embeddings.\\
(4).\quad For each $X\in \mathcal{A}$, there is an exact sequence
\begin{align*}
0 \longrightarrow i_!i^! (X) \longrightarrow X \longrightarrow
j_*j^*(A),
\end{align*}
where the two morphisms are the adjunction morphisms of the
corresponding adjoint pairs. Compare \cite{Par} and \cite[4.1]{FP}.

\vskip 5pt

In the right recollment, one can show that ${\rm Ker}\; j^* ={\rm
Im}\;i_!$ and thus ${\rm Im}\;i_!$ is a Serre subcategory, and via
the functor $j^*$, there exists a natural equivalence of abelian
categories $\mathcal{A}/{{\rm Im}\; i_!}\simeq \mathcal{A}''$.
Moreover since the pair $(i_!, i^!)$ is adjoint, we infer that ${\rm
Im}\;i_!$ is a localizing Serre subcategory. Conversely, given a
localizing Serre subcategory $\mathcal{B}$, one has  a natural exact
sequence $0\longrightarrow \tau(X) \longrightarrow X \longrightarrow
\omega(\pi(X))$ for each $X\in \mathcal{A}$, and note that the two
morphisms involved are the adjunction morphisms of the corresponding
adjoint pairs, therefore we have a right recollement of abelian
categories:
\[\xymatrix@C=40pt{
\mathcal{B} \ar@<+.7ex>[r]^-{\rm inc} & \ar@<+.7ex>[l]^-{\tau}
\mathcal{A} \ar@<+.7ex>[r]^-{\pi} & \mathcal{A}/\mathcal{B}.
\ar@<+.7ex>[l]^-{\omega} }\]
Observe that, in a proper sense, every
right recollement of abelian categories is given in this way.

\subsection{} Verdier's theory of quotient triangulated
categories is analogous to Gabriel's theory. Let $\mathcal{C}$ be a
triangulated category whose shift functor is denoted by $[1]$. Let
$\mathcal{D}\subseteq \mathcal{C}$ be a thick triangulated
subcategory. Set $\Sigma$ to be the class of morphisms $X
\stackrel{f} \longrightarrow Y$ such that it fits in a triangle
$X\stackrel{f} \longrightarrow Y \longrightarrow Z \longrightarrow
X[1]$ with $Z\in \mathcal{D}$. This is a multiplicative system on
$\mathcal{C}$ which is compatible with the triangulation. Recall
that the \emph{quotient category} $\mathcal{C}/\mathcal{D}$ is the
localization of $\mathcal{C}$ with respect to $\Sigma$, and as in
the case of abelian categories, morphisms in
$\mathcal{C}/\mathcal{D}$ are expressed as (right) fractions. Note
that the quotient category $\mathcal{C}/\mathcal{D}$ has a unique
triangulation such that its triangles are precisely those isomorphic
to the image of some triangles in $\mathcal{C}$, and thus the
quotient functor $\pi: \mathcal{C}\longrightarrow
\mathcal{C}/\mathcal{D}$ is a triangle functor. Observe that ${\rm
Ker}\; \pi=\mathcal{D}$. We say that there is a short ``exact
sequence" of triangulated categories
\begin{align*}
\mathcal{D} \stackrel{\rm inc}\longrightarrow \mathcal{C}
\stackrel{\pi}\longrightarrow \mathcal{C}/\mathcal{D}
\end{align*}

\vskip 5pt

 A triangulated analogue of localizing Serre subcategory
is a right admissible subcategory. Recall that a thick triangulated
subcategory $\mathcal{D}\subseteq \mathcal{C}$ is said to be
\emph{right admissible} (or \emph{a Bousfield subcategory}) provided
that the inclusion functor ${\rm inc}: \mathcal{D} \hookrightarrow
\mathcal{C}$ admits a right adjoint (see \cite[Chapter 9]{Ne} and
compare \cite{BK}). The adjoint will be denoted also by $\tau$, and
thanks to Keller (\cite[6.6 and 6.7]{Ke}), the functor $\tau$ is a
triangle functor such that the adjunction morphisms are natural
morphisms of triangle functors. In this case, the quotient functor
$\pi$ also admits a right adjoint functor $\omega:
\mathcal{C}/\mathcal{D} \longrightarrow \mathcal{C}$. Observe that
both $\tau$ and $\omega$ are fully-faithful. Consider the full
subcategory  $\mathcal{D}^\perp=\{X\in \mathcal{C}\; |\; {\rm
Hom}_\mathcal{C}(\mathcal{D}, X)=0\}$, and one observes that this is
a thick triangulated subcategory. It is remarkable that the
essential image ${\rm Im}\; \omega$ of $\omega$ coincides with
$\mathcal{D}^\perp$ and thus the triangulated subcategory
$\mathcal{D}^\perp$ is left admissible. In fact, this gives a
bijection between the class of right admissible subcategories and
the class of left admissible subcategories.

\vskip 5pt

 Another related concept is the right recollement of
triangulated categories. Let $\mathcal{C}, \mathcal{C}',
\mathcal{C}''$ be triangulated categories. Recall that a \emph{right
recollement} of $\mathcal{C}$ with respect to $\mathcal{C}'$ and
$\mathcal{C}''$ is expressed by the following diagram of triangle
functors
\[\xymatrix@C=40pt{
\mathcal{C}' \ar@<+.7ex>[r]^-{i_{!}} & \ar@<+.7ex>[l]^-{i^{!}}
\mathcal{C} \ar@<+.7ex>[r]^-{j^*} & \mathcal{C}''
\ar@<+.7ex>[l]^-{j_*} }\] satisfying the following conditions\\
(1).\quad The pairs $(i_!, i^!)$ and $(j^*, j_*)$ are adjoint
pairs.\\
(2).\quad $j^*i_!=0$. \\
(3).\quad The functors $i_!$ and $j_*$ are full embeddings.\\
(4).\quad For each $X\in \mathcal{C}$, there is a triangle
\begin{align*}
 i_!i^! (X) \longrightarrow X \longrightarrow
j_*j^*(X) \longrightarrow (i_!i^!(X))[1]
\end{align*}
where the left two morphisms are the adjunction morphisms of the
corresponding adjoint pairs. For details, see \cite{BBD, Par} and
\cite[section 3]{Kr}.

\vskip 5pt

In a right recollement, one has ${\rm Ker}\; j^* ={\rm Im}\;i_!$ and
thus ${\rm Im}\;i_!$ is a thick triangulated subcategory, and via
the functor $j^*$, there exists a natural equivalence of
triangulated categories $\mathcal{C}/{{\rm Im}\; i_!}\simeq
\mathcal{C}''$. Moreover since the pair $(i_!, i^!)$ is adjoint, we
infer that ${\rm Im}\;i_!$ is right admissible. Conversely, given a
right admissible subcategory $\mathcal{D}$, one has a natural
triangle $\tau(X) \longrightarrow X \longrightarrow \omega(\pi(X))
\longrightarrow \tau(X)[1]$ for each $X\in \mathcal{C}$, and note
that the left two morphisms are the adjunction morphisms of the
corresponding adjoint pairs, therefore we have a right recollement
of triangulated categories:
\[\xymatrix@C=40pt{
\mathcal{D} \ar@<+.7ex>[r]^-{\rm inc} & \ar@<+.7ex>[l]^-{\tau}
\mathcal{C} \ar@<+.7ex>[r]^-{\pi} & \mathcal{C}/\mathcal{D}.
\ar@<+.7ex>[l]^-{\omega} }\] Note that every right recollement of
triangulated categories is (in a certain sense) given this way.

\vskip 5pt

 Let us end this section by noting that under some conditions a
 localizing Serre subcatgory of an abelian category may
  induce a right admissible subcategory of the derived
 category, or equivalently, a right recollement of derived categories.
 For an abelian category $\mathcal{A}$, denote by
 $D^b(\mathcal{A})$ (resp.  $D^{+}(\mathcal{A})$ ) the derived category of bounded complexes
 (resp. bounded-below complexes) on
 $\mathcal{A}$. A full subcategory $\mathcal{B} \subseteq
 \mathcal{A}$ is said to \emph{have enough $\mathcal{A}$-injective
 objects} provided that for each object $X\in \mathcal{B}$, there is
 a  monomorphism $X\hookrightarrow I$ in $\mathcal{A}$ where $I\in
 \mathcal{B}$ is an injective object in $\mathcal{A}$.

\begin{lem}\label{recollement}
Let $\mathcal{A}$ be an abelian category, $\mathcal{B}$ a localizing
Serre subcategory. Assume that both $\mathcal{A}$ and
$\mathcal{A}/\mathcal{B}$ have enough injective objects, and
$\mathcal{B}$ has enough $\mathcal{A}$-injective objects. Suppose
the corresponding right recollement is given as in {\bf 2.1}. Then
we have a right recollement of derived ctegories:
\[\xymatrix@C=40pt{
D^{+}(\mathcal{B}) \ar@<+.7ex>[r]^-{\rm inc} &
\ar@<+.7ex>[l]^-{R^+\tau} D^{+}(\mathcal{A}) \ar@<+.7ex>[r]^-{\pi} &
D^{+}(\mathcal{A}/\mathcal{B}), \ar@<+.7ex>[l]^-{R^+\omega} }\]
where $R^+\tau$ and $R^+\omega$ are the derived functors, and here
${\rm inc}$ and $\pi$ are applied on complexes term by term.
\end{lem}

\noindent {\bf Proof.}\quad By \cite[Chapter I, Proposition
4.8]{Har}, the natural functor $D^{+}(\mathcal{B}) \longrightarrow
D^{+}_\mathcal{B}(\mathcal{A})$ is a triangle equivalence, where
$D^{+}_\mathcal{B}(A)$ denotes the full triangulated subcategory of
$D^{+}(\mathcal{A})$ consisting of complexes with all the
cohomologies lying in $\mathcal{B}$. By the dual of \cite[Lemma
(1.1)]{CPS}, or the proof of \cite[Corollary 3.3]{Mi}, we infer that
$({\rm inc}, R^+\tau )$ and $(\pi, R^+\omega)$ are adjoint pairs.
Now our recollement follows directly from Miyachi's ``exact
sequence" of derived categories $D^{+}_\mathcal{B}(\mathcal{A})
\stackrel{\rm inc}\longrightarrow D^{+}(\mathcal{A})
\stackrel{\pi}\longrightarrow D^{+}(\mathcal{A}/\mathcal{B})$ (see
\cite[Theorem 3.2]{Mi}).\hfill $\blacksquare$

\section{Gorenstein Categories}
In this section, we will recall the definition of Gorenstein
categories and their basic properties. Before that, we will study
the cohomological dimensions of complexes for later use.

\subsection{} Let us recall the notion of cohomological dimensions in
an abelian category. Let $\mathcal{A}$ be an abelian category and
assume that the bounded derived category $D^b(\mathcal{A})$ is
defined. We will always identify $\mathcal{A}$ as the full
subcategory of $D^b(\mathcal{A})$ consisting of stalk complexes
concentrated at degree $0$. In particular, the $n$-th extension
group of two objects $X, Y$  of $\mathcal{A}$ is defined to be ${\rm
Ext}^n_\mathcal{A}(X, Y):={\rm Hom}_{D^b(\mathcal{A})}(X, Y[n])$ for
any $n\in \mathbb{Z}$, where $[n]$ denotes the $n$-th power of the
shift functor (\cite[p.62]{Har}). Then we have ${\rm
Ext}_\mathcal{A}^n(X, Y)=0$ if $n<0$, and ${\rm
Ext}_\mathcal{A}^0(X, Y)\simeq {\rm Hom}_\mathcal{A}(X, Y)$. Recall
that the \emph{projective dimension} of an object $X$ is defined by
${\rm proj.dim} \; X :={\rm sup}\; \{n\geq  0\; |\; \mbox{ there
exists } Y\in \mathcal{A}, \; {\rm Ext}^n_\mathcal{A}(X, Y)\neq
0\}$. Dually one defines the \emph{injective dimension} of $Y$,
denoted by ${\rm inj.dim}\; Y$. The \emph{global dimension} of the
category $\mathcal{A}$ is defined to be  ${\rm gl.dim}\;
\mathcal{A}:={\rm sup}\; \{{\rm proj.dim}\; X \; |\; X \in
\mathcal{A}\}$. Observe that ${\rm gl.dim}\; \mathcal{A} ={\rm
sup}\; \{{\rm inj.dim}\; Y \; |\; Y \in \mathcal{A}\}$. Also we need
the finitistic dimensions: the \emph{finitistic projective dimension
}of $\mathcal{A}$ is defined to be ${\rm fin.pd}\; \mathcal{A}:={\rm
sup}\; \{{\rm proj.dim}\; X< \infty\; |\; X \in \mathcal{A}\}$;
similarly, we define the \emph{finitistic injective dimension} ${\rm
fin.id}\; \mathcal{A}$. Observe that ${\rm fin.pd}\; \mathcal{A},
{\rm fin.id}\; \mathcal{A}\leq {\rm gl.dim}\; \mathcal{A}$, and if
${\rm gl.dim}\;\mathcal{A} < \infty$, then we have ${\rm fin.pd}\;
\mathcal{A}={\rm gl.dim}\; \mathcal{A}={\rm fin.id}\; \mathcal{A}$.

\vskip 5pt

 For later use, we generalize the above and define the
projective and injective dimensions for bounded complexes. Let
$X^\bullet\in D^b(\mathcal{A})$ be a bounded complex. Define the
\emph{projective dimension} of $X^\bullet$ as
\begin{align*}
{\rm proj.dim}\; X^\bullet:={\rm sup}\; \{n\in \mathbb{Z}\; |\;
\mbox{ there exists } Y \in \mathcal{A}, {\rm
Hom}_{D^b(\mathcal{A})}(X^\bullet, Y[n])\neq 0\}.
\end{align*}
Dually, the \emph{injective dimension} is
\begin{align*}
{\rm inj.dim}\; Y^\bullet:={\rm sup}\; \{n\in \mathbb{Z}\; |\;
\mbox{ there exists } X \in \mathcal{A}, {\rm
Hom}_{D^b(\mathcal{A})}(X[-n], Y^\bullet)\neq 0\}.
\end{align*}
If $X^\bullet=X$ is a  stalk complex concentrated at degree $0$,
then its projective dimension coincides with the one given in the
preceding paragraph. Similar remark holds for $Y^\bullet$. Note that
two quasi-isomorphic complexes necessarily have the same projective
and injective dimensions.

\vskip 5pt

\begin{lem}\label{dimensionofcomplexes}
Let $X^\bullet=(X^n, d^n)_{n\in \mathbb{Z}}$ be a bounded complex.
Then we have ${\rm proj.dim}\; X^\bullet \leq {\rm max}\; \{{\rm
proj.dim}\; X^n -n\; |\; X^n\neq 0\}$ and  ${\rm inj.dim}\;
X^\bullet \leq {\rm max}\; \{{\rm inj.dim}\; X^n +n\; |\; X^n\neq
0\}$.
\end{lem}

\noindent{\bf Proof.}\quad First note two facts: for each complex
$X^\bullet$ and each $n\in \mathbb{Z}$, ${\rm proj.dim}\;
X^\bullet[n]={\rm proj.dim}\; X^\bullet+n$; for each triangle
${X'}^\bullet \longrightarrow {X}^\bullet \longrightarrow
{X''}^\bullet \longrightarrow {X'}^\bullet[1]$ in
$D^b(\mathcal{A})$, applying the cohomological functor ${\rm
Hom}_{D^b(\mathcal{A})} (-, Y[n])$ for any $Y\in \mathcal{A}$ to it,
we obtain that ${\rm Hom}_{D^b(\mathcal{A})} (X^\bullet, Y[n])\neq
0$ implies that ${\rm Hom}_{D^b(\mathcal{A})} (X'^\bullet, Y[n])\neq
0$ or ${\rm Hom}_{D^b(\mathcal{A})} (X''^\bullet, Y[n])\neq 0$,
therefore we have ${\rm proj.dim}\; X^\bullet \leq {\rm max}\;
\{{\rm proj.dim}\; {X'}^\bullet, {\rm proj.dim}\; {X''}^\bullet\}$.
Note that a given bounded complex $X^\bullet$ is an iterated
extension of the stalk complexes $X^n[-n]$, and thus applying the
above two facts repeatedly, we prove the first inequality, and
similarly the second one. \hfill $\blacksquare$ \vskip 5pt

Following \cite[Appendix]{NV},  we denote by $D^b(\mathcal{A})_{\rm
fpd}$ (resp. $D^b(\mathcal{A})_{\rm fid}$) the full subcategory of
$D^b(\mathcal{A})$ consisting of complexes of finite projective
dimension (resp. finite injective dimension). Observe that they are
thick triangulated subcategories. Note that if $\mathcal{A}$ has
enough projective objects (resp. injective objects), the complexes
in $D^b(\mathcal{A})_{\rm fpd}$ (resp. $D^b(\mathcal{A})_{\rm fid}$)
have a characterization by resolutions. For example, we recall that
the abelian category $\mathcal{A}$ is said to \emph{have enough
injective objects} if every object can be embedded into an injective
object. Then the following is well-known (compare \cite[Chapter I,
Proposition 7.6]{Har}).

\begin{lem}\label{finiteinjectivedimension}
Let $\mathcal{A}$ be an abelian category with enough injective
objects, let $X^\bullet=(X^n, d^n) \in D^b(\mathcal{A})$, $n_0 \in
\mathbb{Z}$. The
following are equivalent: \\
(1). \quad ${\rm inj.dim}\; X^\bullet \leq n_0$. \\
(2).\quad For any quasi-isomorphism $X^\bullet \longrightarrow
I^\bullet$ with $I^\bullet=(I^n, \partial^n)$ a bounded-below
complex of injective objects, then ${\rm Ker}\; \partial^{n_0}$ is
injective and
$H^n(I^\bullet)=0$ for $n> n_0$.\\
(3). \quad There exists a quasi-isomorphism $X^\bullet
\longrightarrow I^\bullet$ where $I^\bullet$ is a bounded complex of
injective objects with $I^n=0$ for $n> n_0$. \par
 Consequently, we have a natural equivalence of triangulated categories
 $K^b(\mathcal{I})\simeq D^b(\mathcal{A})_{\rm fid}$, where
 $\mathcal{I}$ is the full subcategory consisting of all the injective
 objects and $K^b(\mathcal{I})$ its bounded homotopy category.
\end{lem}

\noindent {\bf Proof. }\quad For ``$(1) \Longrightarrow (2)$", first
note that $H^n(X^\bullet)=0$ for $n> n_0$. Otherwise, the natural
chain map ${\rm Ker}\; d^n [-n]\longrightarrow X^\bullet$ induces an
epimorphism on the $n$-th cohomology groups, and thus it is not zero
in $D^b(\mathcal{A})$. However this is impossible since ${\rm
inj.dim}\; X^\bullet<n$. Hence via the quasi-isomorphism we have
$H^n(I^\bullet)=0$ for $n> n_0$. Therefore $X^\bullet$ is isomorphic
to the good truncated complex $\tau^{\leq n_0}I^\bullet= \cdots
\longrightarrow I^{n_0-2} \longrightarrow I^{n_0-1} \longrightarrow
{\rm Ker}\; \partial^{n_0}\longrightarrow 0$ in the derived category
$D^b(\mathcal{A})$. Note that we have a natural triangle in
$D^b(\mathcal{A})$
\begin{align*}
(\sigma^{\leq n_0-1}I^\bullet)[-1] \longrightarrow {\rm Ker} \;
\partial^{n_0}[-n_0]\longrightarrow \tau^{\leq n_0}I^\bullet \longrightarrow
\sigma^{\leq n_0-1}I^\bullet
\end{align*}
where $\sigma^{\leq n_0-1} I^\bullet = \cdots \longrightarrow
I^{n_0-2} \longrightarrow I^{n_0-1} \longrightarrow  0 $ is the
brutal truncated complex. Thus since ${\rm inj.dim}\; \tau^{\leq
n_0}I^\bullet={\rm inj.dim}\; X^\bullet \leq n_0 $ and as we will
easily see below in the proof of ``$(3) \Longrightarrow (1)$" that
${\rm inj.dim}\; \sigma^{\leq n_0-1}I^\bullet\leq n_0-1$, and thus
from  the triangle above, we infer that ${\rm inj.dim}\; {\rm Ker}\;
\partial^{n_0}[-n_0]\leq n_0$, and thus since ${\rm inj.dim}\; {\rm Ker}\;
\partial^{n_0}[-n_0]= {\rm inj.dim}\; {\rm Ker}\; \partial^{n_0}
+n_0$, we infer that ${\rm inj.dim}\; {\rm Ker}\;
\partial^{n_0}\leq 0$, and then it necessarily forces  that ${\rm inj.dim}\; {\rm Ker}\;
\partial^{n_0}= 0$, that is, ${\rm Ker}\; \partial^{n_0}$ is
injective.

For ``$(2) \Longrightarrow (3)$", first take a quasi-isomorphism
$X^\bullet \longrightarrow I^\bullet$ where
$I^\bullet=(I^n,\partial^n)$ is a bounded-below complex of injective
objects (\cite[p.42]{Har}). Then by (2), the subcomplex $\tau^{\leq
n_0}I^\bullet = \cdots \longrightarrow  I^{n_0-2}  \longrightarrow
I^{n_0-1} \longrightarrow {\rm Ker}\;
\partial^{n_0} \longrightarrow 0$ is a bounded complex of injective
objects, and moreover the natural chain map $\tau^{\leq
n_0}I^\bullet \longrightarrow I^\bullet$ is split-mono in the
category of complexes. Take its retraction to be $s^\bullet$. Then
the composite $X^\bullet \longrightarrow I^\bullet
\stackrel{s^\bullet}\longrightarrow \tau^{\leq n_0}I^\bullet $
fulfills (3).

 The implication ``$(3) \Longrightarrow (1)$" is
immediate, since we have ${\rm inj.dim}\; X^\bullet={\rm inj.dim}\;
I^\bullet$ and the canonical isomorphism ${\rm
Hom}_{D^b(\mathcal{A})} (X[-n], I^\bullet)\simeq {\rm
Hom}_{K^b(\mathcal{A})} (X[-n], I^\bullet)$ for any $X\in
\mathcal{A}$, here $K^b(\mathcal{A})$ is the bounded homotopy
category of $\mathcal{A}$.\hfill $\blacksquare$

\subsection{}  We will recall the definition of Gorenstein
categories. Recall that an abelian category $\mathcal{A}$ is said to
\emph{satisfy the {\rm (AB4)} condition }if it has arbitrary
(set-indexed) coproducts and coproducts preserve short exact
sequences. An object $T\in \mathcal{A}$ is said to be \emph{a
generator} if the functor ${\rm Hom}_\mathcal{A}(T, -): \mathcal{A}
\longrightarrow {\rm Ab}$ is faithful, where Ab denotes the category
of abelian groups. Observe that if the abelian category
$\mathcal{A}$ has arbitrary coproducts, then an object $T$ is a
generator if and only if any object in $\mathcal{A}$ is a quotient
of a coproduct of copies of $T$. \vskip 5pt

This definition is essentially given by  Enochs, Estrada,
Garc\'{\i}a Rozas (\cite[Definition 2.18]{EEG}). One may note that
there are other different but related notions of Gorenstein
categories, see \cite[section 4]{Be2} and \cite[3.2]{KR}.

\begin{defn}
An abelian category $\mathcal{A}$ is called a \emph{Gorenstein
category},
if the following conditions are satisfied:\\
(G1).\quad $\mathcal{A}$ satisfies (AB4) and has enough injective
objects, $\mathcal{A}$ has a generator $T$ of finite
projective dimension.\\
(G2).\quad An object $X$ has finite projective dimension if and only if it has finite injective dimension. \\
(G3).\quad Every injective object has finite projective dimension
and these dimensions are uniformly bounded.\par Denote by ${\rm
G.dim}\; \mathcal{A} :={\rm max}\; \{{\rm proj.dim} \; I \; | \; I
\mbox{ any injective object}\}$, which is called the
\emph{Gorenstein dimension} of $\mathcal{A}$.
\end{defn}

By the following observation, one deduces that our definition is
equivalent to \cite[Definition 2.18]{EEG} in the case of
Grothendieck categories.

\begin{prop}\label{equivalencesofdefinition}
Let $\mathcal{A}$ be a Gorenstein category. Then we have
\begin{align*}{\rm fin.id}\; \mathcal{A}\leq {\rm
fin.pd}\;\mathcal{A}={\rm G.dim}\; \mathcal{A}.\end{align*}
\end{prop}

\noindent {\bf Proof.}\quad Assume $X\in \mathcal{A}$ and ${\rm
inj.dim}\; X=n_0<\infty$. Then we have an exact sequence
$$0\longrightarrow X \longrightarrow I^0 \longrightarrow \cdots
\longrightarrow I^{n_0-1} \stackrel{\partial}\longrightarrow I^{n_0}
\longrightarrow 0,$$ where each $I^i$ is injective and the
epimorphism $\partial$ is not split, in particular, we have ${\rm
Ext}_\mathcal{A}^1(I^{n_0}, {\rm Ker}\; \partial)\neq 0$. By the
dimension-shift technique in homological algebra, we have ${\rm
Ext}_\mathcal{A}^{n_0}(I^{n_0}, X)\simeq {\rm
Ext}_\mathcal{A}^1(I^{n_0}, {\rm Ker}\;
\partial)\neq 0$. Therefore, $n_0\leq {\rm proj.dim}\; I^{n_0}\leq {\rm
G.dim}\; \mathcal{A}$, and thus ${\rm fin.id}\; \mathcal{A}\leq {\rm
G.dim}\; \mathcal{A}$.

 By (G2) we have ${\rm G.dim}\;\mathcal{A}\leq {\rm fin.pd}\;
 \mathcal{A}$. On the other hand, assume that $X$ has finite
 projective dimension. By (G2), it has finite injective dimension,
 and we take its injective resolution $X\longrightarrow I^\bullet$ as above.
 Hence ${\rm proj.dim}\; X={\rm proj.dim}\; I^\bullet$. Applying Lemma
 \ref{dimensionofcomplexes},
 we get ${\rm proj.dim}\; I^\bullet \leq {\rm max}\; \{{\rm proj.dim}\; I^n-n\; |\; 0\leq n\leq
 n_0\}$, in particular, ${\rm proj.dim}\; I^\bullet\leq {\rm G.dim}\; \mathcal{A}$.
  From this, one infers that ${\rm fin.pd}\; \mathcal{A}\leq {\rm G.dim}\;
 \mathcal{A}$. Then we are done. \hfill $\blacksquare$

\vskip 5pt

 We observe that the condition (G2) can be characterized
in terms of  the derived category (compare \cite[Lemma 1.5
(iii)]{Ha2}).

\begin{prop}\label{derivedGorenstein}
Let $\mathcal{A}$ be an abelian category satisfying (G1). Then the
condition (G2) is equivalent to the condition that
$D^b(\mathcal{A})_{\rm fpd}=D^b(\mathcal{A})_{\rm fid}$.
\end{prop}

\noindent {\bf Proof.}\quad Assume that the condition (G2) holds.
First observe that $D^b(\mathcal{A})_{\rm fid}\subseteq
D^b(\mathcal{A})_{\rm fpd}$. In fact, if $X^\bullet \in
D^b(\mathcal{A})_{\rm fid}$, then by Lemma
\ref{finiteinjectivedimension}, there is a quasi-isomorphism
$X^\bullet \longrightarrow I^\bullet$ with $I^\bullet$ a bounded
complex of injective objects. By (G2) each $I^i$ has finite
projective dimension, and then by Lemma \ref{dimensionofcomplexes},
we obtain that ${\rm proj.dim}\; I^\bullet<\infty$, and thus ${\rm
proj.dim}\; X^\bullet<\infty$. On the other hand, given $X^\bullet
\in D^b(\mathcal{A})_{\rm fpd}$, by the dual of \cite[ Chapter I,
Lemma 4.6(1)]{Har}, we may take a quasi-isomorphism $T^\bullet
\longrightarrow X^\bullet$, where $T^\bullet=(T^n, \delta^n)$ is a
bounded-above complex such that each $T^n$ is a coproduct of copies
of $T$. In particular, ${\rm proj.dim}\; T^n\leq {\rm proj.dim} \;
T< \infty$, and by (G2) ${\rm inj.dim}\; T^n< \infty$. Now we will
follow the idea in the proof of  \cite[Appendix, Lemma A.1]{NV}.
Take $n_0>>0$ such that $H^{-n}(T^\bullet)=0$ for $n\geq n_0$ and
$n_0\geq  {\rm proj.dim}\; X^\bullet$. Consider the good truncated
complex $\tau^{\geq -n_0}T^\bullet:= 0\longrightarrow {\rm Coker}\;
\delta^{-n_0-2}\longrightarrow T^{-n_0} \longrightarrow
T^{-n_0+1}\longrightarrow \cdots $, which is isomorphic to
$X^\bullet$ in $D^b(\mathcal{A})$. Then we have the following
natural triangle
\begin{align*}
\sigma^{\geq -n_0} T^\bullet \longrightarrow \tau^{\geq
-n_0}T^\bullet \longrightarrow {\rm Coker}\;\delta^{-n_0-2}[n_0+1]
\longrightarrow \sigma^{\geq -n_0} T^\bullet[1].
\end{align*}
Since ${\rm proj.dim}\; \tau^{\geq -n_0}T^\bullet\leq n_0$, the
middle morphism in the triangle is zero, and hence the first
morphism is  split-epi, in other words, $\tau^{\geq -n_0}T^\bullet$
is a direct summand of $\sigma^{\geq -n_0} T^\bullet$. As we noted
above that ${\rm inj.dim}\; T^n< \infty$ for each $n$, and then by
Lemma \ref{dimensionofcomplexes}, we have ${\rm inj.dim}\;
\sigma^{\geq -n_0} T^\bullet< \infty$. So we deduce that ${\rm
inj.dim}\; \tau^{\geq -n_0}T^\bullet< \infty$, and since $X^\bullet$
is isomorphic to $\tau^{\geq -n_0}T^\bullet$, we obtain that
$X^\bullet$ also has finite injective dimension.

  Now assume that $D^b(\mathcal{A})_{\rm fpd}=D^b(\mathcal{A})_{\rm
  fid}$. Then (G2) is immediate only if one notes that an object has
  finite projective dimension (resp. finite injective dimension) if and only if it,
  as a stalk complex, lies in $D^b(\mathcal{A})_{\rm
  fpd}$ (resp. $D^b(\mathcal{A})_{\rm
  fid}$). \hfill $\blacksquare$

\vskip 10pt

\section{Quotients of  Gorenstein Categories}

In this section, we will consider the question of  when a quotient
category of a Gorenstein category is still Gorenstein, and we will
prove Theorem A in a general form.

\subsection{} Use the notation in {\bf 2.1}. Assume that the abelian category
$\mathcal{A}$ satisfies (AB4) and the Serre subcategory
$\mathcal{B}$ is closed under arbitrary coproducts, then ones sees
that the corresponding multiplicative system  $\Sigma$ is closed
under coproducts, and furthermore the quotient category
$\mathcal{A}/\mathcal{B}$ has arbitrary coproducts and these
coproducts preserve short exact sequences, that is, the quotient
category $\mathcal{A}/\mathcal{B}$ also satisfies (AB4). Note that
in this case, the functor $\pi$ preserves coproducts.

\vskip 5pt

Our main result in this section is :

\begin{thm}\label{maintheorem}
Let $\mathcal{A}$ be a Gorenstein category with a generator $T$ of
finite projective dimension,  and $\mathcal{B}\subseteq \mathcal{A}$
 a Serre subcategory closed under coproducts,  $\pi: \mathcal{A}\longrightarrow
 \mathcal{A}/\mathcal{B}$ the canonical functor. Assume that the functor $\pi$
 sends injective objects to injective objects and ${\rm proj.dim}\; \pi(T)<
 \infty$. Then the quotient category $\mathcal{A}/\mathcal{B}$ is
 Gorenstein.\par
 Moreover, we have ${\rm G.dim}\; \mathcal{A}/\mathcal{B}\leq {\rm G.dim}\;\mathcal{A}+ {\rm proj.dim}\;
 \pi(T)$.
\end{thm}

\noindent {\bf Proof.}\quad Let us verify the defining conditions
(G1)-(G3) for the quotient category $\mathcal{A}/\mathcal{B}$. For
(G1), as we noted above that $\mathcal{A}/\mathcal{B}$ satisfies
(AB4); and for each object $\pi(X)$, there is a monomorphism $i: X
\longrightarrow I$ in $\mathcal{A}$ with $I$ injective, then by
assumption $\pi(I)$ is injective in which $\pi(X)$ embeds, that is,
$\mathcal{A}/\mathcal{B}$ has enough injectives; note that every
object $X$ is a quotient of a coproduct of copies of $T$, hence
$\pi(X)$ is a quotient of a coproduct of copies of $\pi(T)$ (using
that $\pi$ preserves coproducts), hence $\pi(T)$ is a generator, and
by assumption, it is of finite projective dimension.

For the condition (G2), first note that since $\pi$ is exact and
preserves injectives, $\pi$ sends objects of finite injective
dimension to objects of finite injective dimension. In particular,
if $T'$ is a coproduct of copies of $T$, then by (G2) of
$\mathcal{A}$, it has finite injective dimension, then we infer that
$\pi(T')$ has finite injective dimension. Take $\pi(X) \in
\mathcal{A}/\mathcal{B}$. Assume that ${\rm proj.dim}\; \pi(X)=n_0
<\infty$. Since $\pi(T)$ is a generator, we may take an exact
sequence
\begin{align*}
0\longrightarrow \pi(K) \longrightarrow \pi(T^{-n_0})
\longrightarrow \pi(T^{-n_0+1}) \longrightarrow \cdots
\longrightarrow \pi(T^{0}) \longrightarrow \pi(X) \longrightarrow 0,
\end{align*}
where each $T^{-i}$ is a coproduct of copies of $T$. We have a
natural triangle in $D^b(\mathcal{A}/\mathcal{B})$:
\begin{align*}
\mathcal{T}^\bullet \longrightarrow \pi(X)
\stackrel{\xi}\longrightarrow \pi(K)[n_0+1] \longrightarrow
\mathcal{T}^\bullet[1],
\end{align*}
where $\mathcal{T}^\bullet=0 \longrightarrow \pi(T^{-n_0})
\longrightarrow \pi(T^{-n_0+1}) \longrightarrow \cdots
\longrightarrow \pi(T^{0}) \longrightarrow 0$. Since ${\rm
proj.dim}\; \pi(X)=n_0$, we infer that $\xi$ is zero, that is, the
leftmost morphism is split-epi. Hence $\pi(X)$ is a direct summand
of $\mathcal{T}^\bullet$, hence ${\rm inj.dim}\; \pi(X)\leq {\rm
inj.dim}\; \mathcal{T}^\bullet$. Here again we have used the idea in
the proof of \cite[Appendix, Lemma A.1]{NV}. As we noted above, each
$\pi(T^{-i})$ has finite injective dimension, hence by Lemma
\ref{dimensionofcomplexes}, ${\rm inj.dim}\; \mathcal{T}^\bullet<
\infty$. Therefore $\pi(X)$ has finite injective dimension.

On the other hand, assume that ${\rm inj.dim}\; \pi(X)=n_0<\infty$.
We will show that ${\rm proj.dim}\; \pi(X)< \infty$. Then we have an
exact sequence $0 \longrightarrow \pi(X) \longrightarrow
\mathcal{I}^0 \longrightarrow \cdots \longrightarrow
\mathcal{I}^{n_0-1} \longrightarrow \mathcal{I}^{n_0}
\longrightarrow 0$ with each $\mathcal{I}^i$ injective. By Lemma
\ref{dimensionofcomplexes} it suffices to show that each
$\mathcal{I}^i$, or more generally, any injective object
$\mathcal{I}$ in $\mathcal{A}/\mathcal{B}$ has finite projective
dimension. We claim that $\mathcal{I}$ is direct summand of $\pi(I)$
for some injective object $I$ in $\mathcal{A}$. In fact, take $X'\in
\mathcal{A}$ such that $\pi(X')=\mathcal{I}$, embed $X'$ into an
injective object $I$, thus $\mathcal{I}$ embeds into $\pi(I)$, and
this embedding is necessarily split since $\mathcal{I}$ is
injective, and this proves the claim. So to prove (G2), it suffices
to show that for each injective $I$ in $\mathcal{A}$, ${\rm
proj.dim}\; \pi(I)<\infty$. By (G3) of $\mathcal{A}$, assume that
${\rm proj.dim}\; I = d \leq {\rm G.dim}\; \mathcal{A}$. Then by a
similar argument as above, we have that in $D^b(\mathcal{A})$, $I$
is a direct summand of a complex $T^\bullet= 0 \longrightarrow
T^{-d} \longrightarrow T^{-d+1} \longrightarrow \cdots
\longrightarrow T^0\longrightarrow 0$ with each $T^{-i}$ a coproduct
of copies of $T$. Hence applying $\pi$, we have that in
$D^b(\mathcal{A}/\mathcal{B})$, $\pi(I)$ is a direct summand of
$\pi(T^\bullet)$, and thus ${\rm proj.dim}\; \pi(I) \leq {\rm
proj.dim}\; \pi(T^\bullet)$. Since by assumption $\pi(T)$ has finite
projective dimension, so does any coproducts of its copies, more
precisely, ${\rm proj.dim}\; \pi(T^{-i})\leq {\rm proj.dim}\;
\pi(T)$. Applying Lemma \ref{dimensionofcomplexes}, we obtain that
${\rm proj.dim}\; \pi(T^\bullet)\leq d+{\rm proj.dim}\; \pi(T)$.
Consequently, we have
 ${\rm proj.dim}\; \pi(I)\leq {\rm G.dim}\; \mathcal{A}+{\rm proj.dim}\;
 \pi(T)$. This finishes the proof of (G2) and this also proves  (G3), and
 even more, one infers that ${\rm G.dim}\;\mathcal{A}/\mathcal{B}
 \leq {\rm G.dim}\; \mathcal{A}+{\rm proj.dim}\;
 \pi(T)$. We are done. \hfill $\blacksquare$

\subsection{} In this section, we will apply Theorem \ref{maintheorem} to the
situation where the abelian category $\mathcal{A}$ is Grothendieck
and $\mathcal{B}\subseteq \mathcal{A}$ is a localizing Serre
subcategory.

\vskip 5pt

Before fixing our setup we recall some notions. Recall that a
monomorphism $i: X\longrightarrow Y$ is \emph{essential} if for any
morphism $g: Y \longrightarrow Z$, the composite $g\circ i$ is mono
implies that $g$ is mono; the \emph{injective hull} of $X$ means an
essential monomorphism $i: X\longrightarrow I$ with $I$ injective.
We say a full subcategory $\mathcal{B}$ is \emph{closed under
essential monomorphisms}, if for any essential monomorphism $i:
X\longrightarrow Y$ with $X\in \mathcal{B}$, then $Y$ lies also in
$\mathcal{B}$. Recall that an abelian category $\mathcal{A}$ is said
to be \emph{a Grothendieck category} if it satisfies (AB4), direct
limits in $\mathcal{A}$ are exact, and it has a generator. Note that
in a Grothendieck category, any object has a (unique) injective hull
(\cite[Chapter V]{St}). Assume that $\mathcal{A}$ is Grothendieck
and $\mathcal{B} \subseteq \mathcal{A}$ is a localizing Serre
subcategory, then the quotient category $\mathcal{A}/\mathcal{B}$ is
also Grothendieck (for details, see \cite[Chapter X]{St})

\vskip 10pt

\noindent{\bf Setup 4.2:} \quad  The abelian category $\mathcal{A}$
is Grothendieck, and $\mathcal{B}\subseteq\mathcal{A}$ is a
localizing Serre subcategory which is closed under essential
monomorphisms.

\vskip 10pt

Let us draw some immediate consequences of this setup.

\begin{lem} \label{consequencelemma} Assume Setup 4.2. Then we have  \\
(1).\quad Both the functors ${\rm inc}: \mathcal{B} \longrightarrow \mathcal{A}$  and
$\pi: \mathcal{A} \longrightarrow \mathcal{A}/\mathcal{B}$ preserve injective objects.\\
(2).\quad The subcategory $\mathcal{B}$ has enough
$\mathcal{A}$-injective objects.
\end{lem}

\noindent{\bf Proof.}\quad Let $I$ be any injective object in
$\mathcal{A}$. We first claim that its subobject $\tau(I)$ is
injective. In fact, consider the injective hull $i:
\tau(I)\longrightarrow I'$. Since $i$ is an essential monomorphism,
by the assumption we get that $I'\in \mathcal{B}$. However by the
injectiveness of $I'$, there is a morphism $i': I' \longrightarrow
I$ extending the inclusion of $\tau(I)$ into $I$, which is
necessarily a monomorphism since $i'\circ i$ is mono and $i$ is
essential. Note that $\tau(I)$ is the unique largest subobject of
$I$ belonging to $\mathcal{B}$, hence the image of $i'$ lies in
$\tau(I)$, and consequently we deduce that $i$ is an isomorphism,
and therefore $\tau(I)$ is injective.

 Assume that $I$ is injective as above. Then $\tau(I)$ is also
 injective hence it is direct summand of $I$, say $I\simeq \tau(I)\oplus
 I'$. Then we have that $I'$ is injective and $\tau(I')=0$, and thus $I'$ is $\mathcal{B}$-local,
  by \cite[Lemma A. 2.7]{Ne}, we have a natural isomorphism
  ${\rm Hom}_\mathcal{A}(X, I')\simeq {\rm Hom}_{\mathcal{A}/\mathcal{B}}(\pi(X), \pi(I'))$
  for any $X\in\mathcal{A}$, and note that  any short exact sequence
in $\mathcal{A}/\mathcal{B}$ is isomorphic to the image of some
short exact sequence in $\mathcal{A}$, we infer that the functor
${\rm Hom}_{\mathcal{A}/\mathcal{B}}(-, \pi(I'))$ is exact, and thus
$\pi(I')$ is injective. Note that $\pi(I)\simeq \pi(I')$, we get
that $\pi$ preserves injective objects.

  Let $B\in \mathcal{B}$, take its injective hull in $\mathcal{A}$,
$i: B\longrightarrow I$.  Since $\mathcal{B}$ is closed under
essential monomorphisms, we have that $I$ belongs to $\mathcal{B}$,
and thus this proves (2). Moreover if $B$ is
$\mathcal{B}$-injective, then $i: B\longrightarrow I$ splits and
thus $B$ is a direct summand of $I$, so  we have that $B$ is also
injective in $\mathcal{A}$, in other words, the inclusion functor
${\rm inc}$ preserves injective objects. We are done. \hfill
$\blacksquare$

\vskip 5pt

Recall that the \emph{cohomological dimension} (\cite[p.57]{Har}) of
the functor $\omega: \mathcal{A}/\mathcal{B}\longrightarrow
\mathcal{A}$ is
 defined to be
 \begin{align*}
 {\rm coh.dim}\; \omega := {\rm sup}\; \{n\geq 0\; |\; \mbox{the $n$-th right derived functor }R^n\omega\neq 0
 \}.
 \end{align*}
Similarly one defines the cohomological dimension of $\tau$, denoted
by ${\rm coh.dim}\; \tau$.

\vskip 5pt

 \begin{prop}\label{cohomologicaldimensioncondition}
Assume Setup 4.2. Assume that $\mathcal{A}$ has a generator $T$ of
finite projective dimension. Then ${\rm proj.dim} \; \pi(T)< \infty$
if and only if ${\rm coh.dim}\; \omega<\infty$, if and only if ${\rm
coh.dim}\; \tau <\infty$.
\par
 Moreover, we have ${\rm coh.dim}\; \omega={\rm max}\; \{0, {\rm coh.dim}\;
 \tau-1\}$ and
${\rm coh.dim}\; \omega \leq {\rm proj.dim}\; \pi(T)\leq {\rm
coh.dim}\; \omega+ {\rm proj.dim}\; T$.
 \end{prop}

\noindent {\bf Proof.}\quad By Lemma \ref{consequencelemma}, we can
apply Lemma \ref{recollement} in this case. By the right recollement
in the lemma, we have for each object $X\in \mathcal{A}$, a triangle
in the derived category $R^+\tau(X) \longrightarrow X
\longrightarrow R^+\omega( \pi(X)) \longrightarrow R^+\tau(X)[1]$,
and taking their cohomologies, we infer that $R^{i+1}\tau(X)\simeq
R^i\omega(\pi(X))$ for each $i\geq 1$. Thus it follows that  ${\rm
coh.dim}\; \omega<\infty$ if and only if ${\rm coh.dim}\; \tau
<\infty$, and moreover, we have ${\rm coh.dim}\; \omega={\rm max}\;
\{0, {\rm coh.dim}\;
 \tau-1\}$ (compare \cite[p.521, line 13]{NV}).

 To continue the proof, let us recall two  facts: let
$I^\bullet=(I^n, d^n)_{n\in \mathbb{Z}}$ be a complex of injective
objects in $\mathcal{A}$. The first fact is, for any fixed $n$, if
$H^n(I^\bullet)\neq 0$, then $H^n({\rm Hom_\mathcal{A}}(T,
I^\bullet))\neq 0$. In fact, in this case, the natural inclusion
${\rm Im}\; d^{n-1}\longrightarrow {\rm Ker}\; d^n$ is proper, thus
since $T$ is a generator, there exists $f: T\longrightarrow {\rm
Ker}\; d^n$ which does not factor through ${\rm Im}\; d^{n-1}$. Then
one observes that this $f$ gives a nonzero element in the cohomology
group $H^n({\rm Hom_\mathcal{A}}(T, I^\bullet))$. The second is
that, if there exists $n_0$ such that $H^n(I^\bullet)=0$ whenever
$n\geq n_0$, then $H^n({\rm Hom_\mathcal{A}}(T, I^\bullet))=0$ for
all $n\geq n_0+{\rm proj.dim}\; T$. Note that the brutal truncated
complex $\sigma^{\geq n_0-1}I^\bullet$ is a shifted version of an
injective resolution of ${\rm Ker}\;d^{n_0-1}$, and thus we have
$H^n({\rm Hom_\mathcal{A}}(T, I^\bullet)) \simeq {\rm
Ext}^{n-n_0+1}_\mathcal{A}(T, {\rm Ker}\;d^{n_0-1})$, and thus the
fact follows.

Take an arbitrary object $\mathcal{M}\in \mathcal{A}/\mathcal{B}$,
and take its injective resolution $\mathcal{M}\longrightarrow
\mathcal{I}^\bullet$. Then we have natural isomorphisms, for all
$n\geq 0$,
\begin{align*}
{\rm Ext}_{\mathcal{A}/\mathcal{B}}^n(\pi(T), \mathcal{M}) =H^n({\rm
Hom}_{\mathcal{A}/\mathcal{B}}(\pi(T), \mathcal{I^\bullet})) \simeq
H^n({\rm Hom}_\mathcal{A}(T, \omega(\mathcal{I^\bullet}))).
\end{align*}
Set $I^\bullet=\omega(\mathcal{I}^\bullet)$, which is a
bounded-below complex of injective objects by Lemma
\ref{consequencelemma} (1). Note that
$H^i(I^\bullet)=R^i\omega(\mathcal{M})$. Then one deduces the result
from the above recalled two facts immediately. \hfill $\blacksquare$

\vskip 10pt

The following is a direct consequence of Theorem \ref{maintheorem}.

\begin{cor}\label{improtantcorollary}
Assume Setup 4.2. Suppose furthur that the category $\mathcal{A}$ is
Gorenstein and the torsion functor $\tau$ has finite cohomological
dimension. Then the quotient category $\mathcal{A}/\mathcal{B}$
Gorenstein.
\end{cor}

\noindent {\bf Proof.}\quad Just note that localizing subcategory is
always closed under coproducts. Now by Lemma \ref{consequencelemma}
(1) and Proposition \ref{cohomologicaldimensioncondition}, we may
apply Theorem \ref{maintheorem}, and thus we are done. \hfill
$\blacksquare$

\section{A Right Recollement for Gorenstein-Injective Objects}

In this section, we will give a right recollement of the stable
category of Gorenstein-injective objects under the same conditions
of Corollary \ref{improtantcorollary}, proving Theorem B in a
general form.

\subsection{} We will recall the concept of Gorenstein-injective
objects. Let $\mathcal{A}$ be an abelian category with enough
injective objects for this moment. Recall that a complex $I^\bullet$
of injective objects is said to be \emph{totally-acyclic}, if it is
exact (= acyclic) and for each injective object $Q$, the Hom complex
${\rm Hom}_\mathcal{A}(Q, I^\bullet)$ is exact. An object $G\in
\mathcal{A}$ is said to be \emph{Gorenstein-injective} provided
there exists a totally-acyclic complex $I^\bullet$  such that
$G=Z^0(I^\bullet)$ is the $0$-th cocycle, and in this case,
$I^\bullet$ is said to be a \emph{complete resolution }of $G$.
Denote by ${\rm GInj}(\mathcal{A})$ the full subcategory of
Gorenstein-injective objects (\cite[Chapter 10]{EJ} and
\cite[section 7]{Kr}). Note that this is an additive subcategory and
injective objects are Gorenstein-injective. Observe that for each
Gorenstein-injective object $G$, we have ${\rm Ext}_\mathcal{A}^i(Q,
G)=0$ for any injective object $Q$ and $i\geq 1$. (In fact, view the
bruntal truncated complex $\sigma^{\geq 0}I^\bullet$ as an injective
resolution of $G$, and then we see that ${\rm Ext}_\mathcal{A}^i(Q,
G)$ is just the $i$-th cohomology group of the Hom complex ${\rm
Hom}_\mathcal{A}(Q, I^\bullet)$, which is zero by the assumption.)

\vskip 5pt

 We collect the basic properties of the category ${\rm
GInj}(\mathcal{A})$.

\begin{lem}\label{basicproperty}
With the notation as above. We have \\
(1).\quad The full subcategory ${\rm GInj}(\mathcal{A})$ is closed
under cokernels of monomorphisms, extensions and taking direct
summands.\\
(2).\quad Endow the exact structure on ${\rm GInj}(\mathcal{A})$ by
short exact sequences in $\mathcal{A}$. Then ${\rm
GInj}(\mathcal{A})$ is a Frobenius category in the sense of
\cite{Ha1} (and \cite[p.381]{Ke3}), whose relative
injective-projective objects are precisely the injective objects in
$\mathcal{A}$.
\end{lem}

\noindent {\bf Proof.}\quad The first statement is proved in
\cite[Theorem 10.1.4]{EJ} and also by the dual of \cite[Proposition
5.1]{AR}. Since ${\rm GInj}(\mathcal{A})$ is closed under
extensions, then it becomes an exact category in the sense of
Quillen by identifying conflations (= admissible exact sequences)
with short exact sequences in $\mathcal{A}$ with terms inside ${\rm
GInj}(\mathcal{A})$, see \cite[Appendix A]{Ke3}. Clearly injective
objects are relative injective and thus ${\rm GInj}(\mathcal{A})$
has enough relative injectives.  Note that ${\rm
Ext}_\mathcal{A}^1(Q, G)=0$ for injective objects $Q$ and $G\in {\rm
GInj}(\mathcal{A})$, and this proves that the functor ${\rm
Hom}_\mathcal{A}(Q, -)$ is exact on all the conflations, i.e., $Q$
is a relative projective. Observe that for any $G\in{\rm
GInj}(\mathcal{A})$, there exists an exact sequence
$0\longrightarrow G' \longrightarrow Q\longrightarrow
G\longrightarrow 0$ with $G'$ Gorenstein-injective and $Q$ injective
in $\mathcal{A}$, and since as we just saw, $Q$ is relative
projective, therefore ${\rm GInj}(\mathcal{A})$ has enough relative
projectives. Now it is not hard to see (2). \hfill $\blacksquare$
\vskip 5pt

\vskip 5pt

The following important result could be derived from \cite[Theorem
2.24]{EEG}.

\begin{prop}\label{keryproposition} Let $\mathcal{A}$ be a Grothendick category which is
Gorenstein, $M\in \mathcal{A}$ an object.  Then the following are
equivalent:\\
(1).\quad The object $M$ is Gorenstein-injective.\\
(2).\quad There exists a long exact sequence $\cdots \longrightarrow
I^{-n}\longrightarrow \cdots \longrightarrow I^{-1} \longrightarrow
I^0\longrightarrow M\longrightarrow 0$ with each $I^{-n}$ injective. \\
(3).\quad For each injective $Q$ and $i\geq 1$, ${\rm
Ext}_\mathcal{A}^i(Q, M)=0$.\\
(4).\quad For each object $L$ of finite injective dimension, ${\rm
Ext}_\mathcal{A}^1(L, M)=0$.\\
\end{prop}

\noindent {\bf Proof.}\quad The implication ``$(1)\Longrightarrow
(2)$" is trivial from the definition. To see ``$(2)\Longrightarrow
(3)$", we apply the dimension-shift technique in homological algebra
and then we get, for any $i\geq 1$ and $k\geq 0$, ${\rm
Ext}_\mathcal{A}^i(Q, M)\simeq {\rm Ext}_{\mathcal{A}}^{i+k+1}(Q,
Z^k)$, where $Z^k$ is the $-k$-th cocycle of the long exact complex.
Note that by (G2), $Q$ has finite projective dimension, and we infer
that ${\rm Ext}_\mathcal{A}^i(Q, M)=0$. ``$(3)\Longrightarrow (4)$"
could also be shown by  the dimension-shift technique. While
``$(4)\Longrightarrow (1)$" is exactly stated in \cite[Theorem
2.24]{EEG}. \hfill $\blacksquare$

\subsection{} Let $\mathcal{A}$ be an abelian category with enough
injectives and let ${\rm GInj}(\mathcal{A})$ be its subcategory of
Gorenstein-injective objects. Consider the stable category
$\underline{{\rm GInj}}(\mathcal{A})$ modulo injective objects: the
objects are the same as in ${\rm GInj}(\mathcal{A})$, and the
morphism space is the quotient of the corresponding one modulo those
factoring through injective objects. For each morphism $f: X
\longrightarrow Y$, we write the corresponding morphism in
$\underline{{\rm GInj}}(\mathcal{A})$ as $\underline{f}: X
\longrightarrow Y$. Since ${{\rm GInj}}(\mathcal{A})$ is a Frobenius
category, then by \cite[Chapter 1, section 2]{Ha1}, the stable
category $\underline{{\rm GInj}}(\mathcal{A})$ has a canonical
triangulated structure as follows:  for each $X\in {\rm
GInj}(\mathcal{A})$, fix a short exact sequence $0\longrightarrow X
\stackrel{i_X}\longrightarrow I(X) \stackrel{d_X}\longrightarrow
S(X) \longrightarrow 0$ such that $I(X)$ is injective and thus
$S(X)$ belongs to ${\rm GInj}(\mathcal{A})$, and thus we have the
induced functor $S: \underline{{\rm
GInj}}(\mathcal{A})\longrightarrow \underline{{\rm
GInj}}(\mathcal{A})$ which is an auto-equivalence and will be the
shift functor; triangles in $\underline{{\rm GInj}}(\mathcal{A})$
are induced by conflations, i.e., the short exact sequences in ${\rm
GInj}(\mathcal{A})$, more precisely, given a short exact sequence
$0\longrightarrow X \stackrel{u} \longrightarrow Y
\stackrel{v}\longrightarrow Z \longrightarrow 0$, we have a
commutative diagram  in $\mathcal{A}$ (using the relative
injectivity of $I(X)$)
\[\xymatrix{
0 \ar[r] & X \ar[r]^-u \ar@{=}[d] & Y \ar[r]^-v \ar@{.>}[d] & Z
\ar[r] \ar@{.>}[d]^-{w} & 0 \\
0 \ar[r] & X \ar[r]^-{i_X} & I(X) \ar[r]^-{d_X} &S(X)\ar[r] & 0. }\]
Then $X \stackrel{\underline{u}} \longrightarrow Y
\stackrel{\underline{v}} \longrightarrow Z
\stackrel{-\underline{w}}\longrightarrow S(X)$ is a triangle, and
all triangles arise in this way (\cite[Lemma 1.2]{CZ}).

\vskip 5pt

 Our main observation is that under the conditions of
Corollary 4.4, we have a natural right recollement of triangulated
categories relating the Gorenstein-injective objects of the abelian
category and of its quotient category.

\begin{thm}\label{maintheoremII} Assume Setup 4.2. Assume that $\mathcal{A}$ is Gorenstein and the torsion
functor $\tau$ has finite cohomological dimension. Then we have a
right recollement of triangulated categories:
\[\xymatrix@C=40pt{
\mathcal{\underline{\rm GInj}(\mathcal{B})} \ar@<+.7ex>[r]^-{\rm
inc} & \ar@<+.7ex>[l]^-{\underline{\tau}} \mathcal{\underline{\rm
GInj}(\mathcal{A})} \ar@<+.7ex>[r]^-{\underline{\pi}} &
\ar@<+.7ex>[l]^-{\underline{\omega}} \underline{\rm
GInj}(\mathcal{A}/\mathcal{B}), }\] where the functors are induced
from the ones between the abelian categories.
\end{thm}

\noindent {\bf Proof.}\quad We will first show that the four
functors involved are well-defined.  This will take several steps.

\vskip 3pt

 We  claim that for each $i\geq 1$ and $X\in {\rm
GInj}(\mathcal{A})$, $R^i\tau(X)=0$. Consequently, $\tau$ preserves
 short exact sequences in ${\rm GInj}(\mathcal{A})$. In fact, since
$X$ is Gorenstein-injective, we have an exact sequence $\cdots
\longrightarrow I^{-n}\stackrel{d^{-n}}\longrightarrow \cdots
\longrightarrow I^{-1} \stackrel{d^{-1}}\longrightarrow I^0
\longrightarrow X\longrightarrow 0$ with each $I^{-n}$ injective.
Thus by the dimension-shift technique, we get $R^i\tau(X)\simeq
R^{i+k}\tau({\rm Im}\; d^{-k})$ for any $k\geq 1$. However $\tau$
has finite cohomological dimension, thus we are done with
$R^i\tau(X)=0$, $i\geq 1$. The other statement follows directly from
the long exact sequence associated with a given short exact sequence
and the derived functors $R^i\tau$.

Next we claim that ${\rm GInj}(\mathcal{B})\subseteq {\rm
GInj}(\mathcal{A})\cap \mathcal{B}$. To see this, let $X\in {\rm
GInj}(\mathcal{B})$, and then from the definition of
Gorenstein-injective objects, we have an exact sequence  $\cdots
\longrightarrow I^{-n}\longrightarrow \cdots \longrightarrow I^{-1}
\longrightarrow I^0\longrightarrow X\longrightarrow 0$ with each
$I^{-n}$ injective in $\mathcal{B}$, however as we observed in Lemma
\ref{consequencelemma} (1), each $I^{-n}$ is also injective in
$\mathcal{A}$, and since $\mathcal{A}$ is Gorenstein and by
Proposition \ref{keryproposition} (2), $X$ belongs to ${\rm
GInj}(\mathcal{A})$.

We claim that if  $X\in {\rm GInj}(\mathcal{A})$, then $\tau(X)\in
{\rm GInj}(\mathcal{B})$. Consequently, combining above we infer
that ${\rm GInj}(\mathcal{B})= {\rm GInj}(\mathcal{A})\cap
\mathcal{B}$. Take the complete resolution $I^\bullet=(I^n, d^n)$
with $Z^0(I^\bullet)=X$. Break $I^\bullet$ into short exact
sequences in ${\rm GInj}(\mathcal{A})$ and note that by the first
claim $\tau$ preserves such short exact sequences, we obtain that
the complex  $\tau(I^\bullet)$ is exact and
$Z^0(\tau(I^\bullet))=\tau(X)$. Note that $\tau$ preserves injective
objects (which is a consequence of a general fact that the right
adjoint functor of an exact functor preserves injectives), hence
$\tau(I^\bullet)$ is an exact sequence of injective objects in
$\mathcal{B}$. Given any injective $Q$ in $\mathcal{B}$, then by
Lemma \ref{consequencelemma} (1), $Q$ is also injective in
$\mathcal{A}$. Note that by adjoint we have an isomorphism of Hom
complexes ${\rm Hom}_\mathcal{B}(Q, \tau(I^\bullet))\simeq {\rm
Hom}_\mathcal{A}(Q, I^\bullet)$, and since $I^\bullet$ is
totally-acyclic, and hence ${\rm Hom}_\mathcal{A}(Q, I^\bullet)$ is
cyclic , and so ${\rm Hom}_\mathcal{B}(Q, \tau(I^\bullet))$ is also
acyclic, in other words, $\tau(I^\bullet)$ is a totally-acyclic
complex in $\mathcal{B}$ and note that
$Z^0(\tau(I^\bullet))=\tau(X)$,  therefore $\tau(X)\in {\rm
GInj}(\mathcal{B})$.

 Let us summarize what we have shown: we have two well-defined
 functors ${\rm inc}: {\rm GInj}(\mathcal{B}) \longrightarrow {\rm GInj}(\mathcal{A})$
and $\tau: {\rm GInj}(\mathcal{A}) \longrightarrow {\rm
GInj}(\mathcal{B})$, both of which preserve relative
injective-projective objects and short exact sequences. Therefore we
have two induced functors on the stable categories and they are
triangle functors by \cite[p.23, Lemma]{Ha1}. Obviously the obtained
pair $({\rm inc}, \underline{\tau})$ is adjoint and the functor
${\rm inc}$ is fully-faithful. \vskip 3pt

 To continue, we claim that
if $X\in {\rm GInj}(\mathcal{A})$, then $\pi(X)\in {\rm
GInj}(\mathcal{A}/\mathcal{B})$. Take the complete resolution
$I^\bullet$ for $X$, then by Lemma \ref{consequencelemma}(1)
$\pi(I^\bullet)$ is an exact sequence of injective objects and its
$0$-th cocycle is $\pi(X)$. By Corollary \ref{improtantcorollary},
$\mathcal{A}/\mathcal{B}$ is Gorenstein. Applying Proposition
\ref{keryproposition}(2), we get that $\pi(X)$ is
Gorenstein-injective.

We claim that for each $i\geq 1$ and $\mathcal{M} \in {\rm
GInj}(\mathcal{A}/\mathcal{B})$, $R^i\omega(\mathcal{M})=0$.
Consequently, $\omega $ preserves  short exact sequences in ${\rm
GInj}(\mathcal{A}/\mathcal{B})$. Note that by Proposition
\ref{cohomologicaldimensioncondition}, the functor $\omega$ has also
finite cohomological dimension, and the proof is similar as the
first claim, and left to the reader.

We claim that if $\mathcal{M}\in {\rm
GInj}(\mathcal{A}/\mathcal{B})$, then $\omega(\mathcal{M})\in {\rm
GInj}(\mathcal{A})$. Note that by Lemma \ref{consequencelemma} (1),
$\pi$ preserves injective objects, and using the last claim, we can
prove this by an argument similar to that in the third claim. And so
the proof is left to the reader.

Hence we have obtained two well-defined functors: $\pi: {\rm
GInj}(\mathcal{A}) \longrightarrow {\rm
GInj}(\mathcal{A}/\mathcal{B})$ and $\omega: {\rm
GInj}(\mathcal{A}/\mathcal{B}) \longrightarrow {\rm
GInj}(\mathcal{A})$, both of which preserve relative
injective-projective objects and short exact sequences. Therefore we
get two induced triangle functors $\underline{\pi}$ and
$\underline{\omega}$ on the stable categories, and the pair
$(\underline{\pi}, \underline{\omega})$ is adjoint. Note that the
functor $\underline{\omega}$ is obviously fully-faithful and the
composite $\underline{\pi}\; {\rm inc}=0$.

\vskip 3pt

To end the proof, it suffices to check the last defining condition
of a right recollement. In fact, from the right recollement in Lemma
\ref{recollement}, we get for each $X\in \mathcal{A}$, a triangle in
the derived category $R^+\tau (X) \longrightarrow X \longrightarrow
R^+\omega(\pi(X)) \longrightarrow R^+\tau(X)[1]$. From the triangle
we get a long exact sequence of cohomological groups
\begin{align*}0
\longrightarrow \tau(X) \longrightarrow X \longrightarrow
\omega(\pi(X)) \longrightarrow R^1\tau(X) \longrightarrow \cdots
\end{align*}
Take $X\in {\rm GInj}(\mathcal{A})$ and by the first claim, we get
an exact sequence in $\mathcal{A}$, even in ${\rm
GInj}(\mathcal{A})$, $0\longrightarrow \tau(X) \longrightarrow X
\longrightarrow \omega(\pi(X)) \longrightarrow 0$, and note that the
morphisms involved are the adjunction morphisms. Short exact
sequences induce triangles, and thus we get the required triangle
for the last defining condition of a right recollement. Then we are
done. \hfill $\blacksquare$

\vskip 20pt

\noindent {\bf Acknowledgement:} \quad The author would like to
thank Prof. Edgar E. Enochs very much for many helpful private
communications during the preparation of this work. Also the author
would to thank Prof. Michel Van den Bergh for his help on the
remarks  concerning local cohomology.

\bibliography{}

\end{document}